\documentclass[12pt,english]{article} 

\small \normalsize
             
\begin{document}
\begin{center}
{\Large\bf A New Look at the Visual Performance of\\[0.3cm]
Nonparametric
Hazard Rate Estimators}\\[0.5cm]
O.~Gefeller$^1$, N.L.~Hjort$^2$ \\[0.3cm]
{\small $^1$ Department of Medical Statistics, University of G\"ottingen,\\ 
Humboldtallee 32, D--37073 G\"ottingen, Germany\\[0.2cm]
$^2$ Department of Mathematics and Statistics, University of Oslo,\\
P.B. 1053 Blindern, N--0316 Oslo, Norway}
\end{center}

\centerline{\bf June 1997} 

\medskip


{\small{\bf Abstract:}            
Nonparametric curve estimation by kernel methods has attracted widespread
interest in theoretical and applied statistics.
One area of conflict between
theory and application relates to the evaluation of the performance
of the estimators. 
Recently, Marron and
Tsybakov (1995) proposed {\it visual error criteria} for addressing this
issue of controversy in density estimation. Their core idea 
consists in using integrated alternatives to the Hausdorff
distance for measuring the closeness of two sets based on
the Euclidean distance.
In this paper, we transfer these ideas to hazard rate estimation from censored
data. We are able to derive similar results that
help to understand when the application
of the new criteria will lead to answers that differ from 
those given by the conventional
approach.}

\section{Introduction}

In various areas of applied 
and theoretical statistics nonparametric 
curve estimation by kernel methods has attracted widespread
interest during the last two decades.
Without imposing any distributional assumptions on the observed data,
structural information, for example, 
on the underlying density functions, some
interesting functionals of the density (like the hazard rate) or
regression curves, can be obtained by ``smoothing" the empirical mass
of the observations to some neighbouring environment around the observed
data points. Starting with the pioneering work by Rosenblatt (1956) and
Parzen (1962), nowadays a vast literature on the properties of kernel
methods can be found (for recent textbooks see Wand and
Jones (1995), H{\"a}rdle (1991)).
Important progress has been made recently
in a variety of issues (for example, bandwidth
selection, boundary behaviour, software implementation), mostly originating
from research in density estimation, the simplest situation. However,
it has been criticised repeatedly that there is an obvious gap between
theory and application with respect to the evaluation of the performance
of the nonparametric estimators. The estimated curves that are judged
to be ``good'' from theoretical reasoning (based on their distance to the
true curve in terms of the classical norms on function spaces, mostly
$L_2$) do not satisfy applied statisticians as their ``graphical fit''
can be substantially poorer than that of other estimators which are
more appealing from a graphical notion but exhibit larger deviations from
the true curve in $L_2$-norm. To address this issue of conflict, Marron and
Tsybakov (1995) proposed new {\it visual error criteria} for assessing the
performance of nonparametric density estimators. The core idea of this
concept consists of using integrated alternatives to the Hausdorff
distance for measuring the closeness of two sets of points based on
the usual Euclidean distance.

The aim of this paper is twofold: first, we review this interesting
approach and describe the basic properties of the visual error criteria
in section 2 and 3.
Second, in the following sections 4 and 5 of this
paper, we apply the criteria to nonparametric hazard rate estimation
from censored data via kernel estimators. The evolving area 
of hazard rate estimation is
a promising field for these new ideas, because, due to the
application-orientated interpretation of hazard rates in practical
survival analysis, qualitative aspects of smoothing performance should
be preferred to the exclusive consideration of $L_2$-fit. As a
first step into this direction, we derive similar results
as in the density context showing that a suitable
asymptotic bridge can be built between expected
visual error criteria and conventional integrated mean squared errors ($MISE$),
incorporating some weight function depending on the derivative of
the hazard rate
In the final section we discuss advantages and disadvantages of this
approach and point to future work needed to broaden the
understanding of the concept.

\section{Conventional \boldmath$L_p$-norm Criteria\unboldmath}
The conventional approach for measuring the distance between some
``true" curve $h$ and its estimate $\widehat{h}$ relies on the application
of classical mathematical norms on function spaces. The most popular
representative of this approach
is the $L_2$-norm. The distance between two functions in
terms of the $L_2$-norm is then defined as
\[
L_2 (\widehat{h}, h) := \int \left[ \widehat{h} (x) - h(x)\right]^2 dx\, .
\]
Alternatively, other variants of the broad class of $L_p$-norms such as,
for example,
$L_1 (\widehat{h}, h) := \int | \widehat{h} (x) - h(x)| dx$ and
$L_{\infty} (\widehat{h}, h) := \sup_x |\widehat{h} (x) - h(x)|$ 
have also been proposed, but are less often used.
All measurements of distance based on the various $L_p$-norms have one
fundamental property in common: they focus completely on the {\it vertical}
discrepancy between the curves at all points $x$ belonging to the
support of the functions (typically in cases of lifetime data, the support
is some subset of ${\cal R}_+$). 

A striking graphical example that the
concentration on vertical discrepancy can contradict the visual notion
of distance between curves has been provided by Kooperberg and Stone (1991).
Their example shows a ``true" bimodal curve that is approximated by two
very different estimates. The key element of this example is
that the visually more appealing estimate that recovers the bimodal
structure but with some imprecision in the location of the second peak
is clearly inferior in terms of all $L_p$-norms to the other estimate
that only recovers the first peak and ``smoothes away" the lower second
peak. The reason for this contradiction between graphical impression and
mathematical $L_p$-norm analysis is obvious: as the distances between the
curves are only assessed vertically, no $L_p$-norm criterion realises that
the first estimate is ``structurally correct", although quite incorrect
with respect to localising the second peak, and thus penalises the estimate
in both regions of vertical discrepancy, the region of the true second peak
and that of the estimated second peak, whereas the ``oversmoothed" estimate
is penalised only once, in the region of the true peak. Any
remedy to this problem has to drop the strict reliance on vertical
distances and has to offer another way of conceptualising discrepancies
between curves.

\section{Visual Error Criteria}
One such alternative approach of measuring the discrepancy between curves
has been developed by Marron and Tsybakov (1995) who termed their measures
{\it visual error criteria}. The starting point for the derivation of these
criteria consists in treating curves not as pointwise collections of the
values of functions $\widehat{h}$ and $h$ of a single variable 
at fixed $x \in {\cal R}_+$ but instead
as sets in ${\cal R}^2$ that are defined by their corresponding ``graphs".
The graph of some function $h$ is defined as the set
$G_h := \{ (x, y), x \in {\cal R_+}, y = h(x)\} \subset {\cal R}^2$.
Some planar distance between the sets $G_h$
and $G_{\widehat{h}}$ can now replace the conventional vertical distance
between $h$ and $\widehat{h}$ for fixed $x \in {\cal R}_+$. The basis
of the definition of a planar distance is the notion of distance from a
point to a set, defined as
\[
d ((x,y), G_h) := \inf_{(x^* , y^*) \in G_h} \| (x, y) - (x^* , y^*)\|_{_2}
\]
and giving the shortest distance from some fixed point $(x,y)$ to any
element in the set $G_h$, where $\|\cdot\|_2$ 
denotes the usual Euclidean distance.
By allowing $(x,y)$ to vary and take on all elements in the set $G_{\widehat{h}}$,
a set of distances between between $G_h$ and $G_{\widehat{h}}$ is defined,
formally introduced as ${\cal M} (G_h , G_{\widehat{h}}) := 
\{ d((x,y), G_{\widehat{h}}): (x,y) \in G_h\}$. A variety of different
ways to summarise 
the complete information in ${\cal M} (G_h , G_{\widehat{h}})$
when defining a real-valued one-dimensional distance are now conceivable.
In analogy to the $L_2$-norm criterion based on quadratic vertical
distances we first consider the two versions of the quadratic visual
error criteria:
\begin{eqnarray*}
VE_2 (h \longrightarrow \widehat{h})& :=&  \left[\int d((x,h(x)), 
G_{\widehat{h}})^2 dx\right]^{1/2}\\
VE_2 (\widehat{h} \longrightarrow h)& :=& \left[\int d((x,\widehat{h}(x)), 
G_h)^2 dx\right]^{1/2}
\end{eqnarray*}
As the planar distances are asymmetrical in nature, both versions have to be
distinguished since they are not identical in general. However, a
symmetrised version of $VE_2$ can also be defined by simply averaging
the two versions, in a Pythagorean way, as follows
\[
SE_2 (\widehat{h}, h ) := \left[ VE_2 (h \longrightarrow \widehat{h})^2
+ VE_2 (\widehat{h} \longrightarrow h)^2 \right]^{1/2}
\]
As for the $L_p$-norm criteria, similar asymmetrical and symmetrical
visual error criteria can also be defined for the non-quadratic distances
(by simply integrating over $d((x,h(x)), G_{\widehat{h}})$ and
$d((x, \widehat{h}(x)), G_h)$, respectively, yielding $VE_1 (h\longrightarrow
\widehat{h})$ and $VE_1 (\widehat{h} \longrightarrow h)$, respectively)
and maximal distance (by considering $\sup ({\cal M}(G_h , 
G_{\widehat{h}}))$ and $\sup ({\cal M}(G_{\widehat{h}}, G_h))$,
respectively, yielding $VE_{\infty}(h \longrightarrow \widehat{h})$ and
$VE_{\infty} (\widehat{h} \longrightarrow h)$, respectively). The
symmetrical version of the two variants of $VE_{\infty}$ given by
$SE_{\infty} (\widehat{h}, h) :=$
\mbox{$\max \{VE_{\infty} (h \longrightarrow \widehat{h}),$}
$VE_{\infty} (\widehat{h} \longrightarrow h)\}$ is also known for a long
time as the Hausdorff distance. Thus, all other visual error criteria
described in this section can be viewed as integrated alternatives to
the Hausdorff distance. 

It has to be kept in mind that the other
symmetrised versions ($SE_2$ and $SE_1$) are --- in a mathematical sense ---
not ``distances" on the corresponding function space as they do not satisfy
the triangle inequality (see Marron and Tsybakov (1995) for an illustrative
counterexample). For the practical application of these criteria,
lacking this mathematical property does not seem to be a serious drawback.

\section{Kernel Estimation of the Hazard Rate}
The evolving area of hazard rate estimation from censored data comprises
a promising field for the application of the visual error criteria as
qualitative aspects of smoothing performance are of primary interest here.
In this section the necessary background about the statistical setting
and some known asymptotic results concerning kernel estimators of the
hazard rate are briefly summarised. More detailed reviews on this subject
can be found in Gefeller and Michels (1992), Andersen et
al. (1993) and Hjort (1996).

Suppose $T_1, \ldots , T_n$ refer to i.i.d. nonnegative failure times
with distribution function $F$ and density function $f$, and
$C_1, \ldots , C_n$ denote i.i.d. nonnegative censoring times with
distribution function $G$ and density function $g$. Assume further that
failure times $T_i$ and censoring times $C_i$ are independent for
all $i=1, \ldots ,n$. Under this setting of the so-called
{\it random censorship model}, which is the simplest and most important
special case of models for censored failure time data in the framework
of counting process models, the observed data consist of the bivariate
sample $(X_1, \delta_1), \ldots , (X_n, \delta_n)$, where
$X_i := \min (T_i, C_i)$ and $\delta_i := I\{T_i \leq C_i\}, i=1, \ldots ,n$.
The censoring indicator $\delta_i$ provides the information whether the
observed $X_i$ refers to a true failure time $(\delta_i = 1)$ or to a
censoring time $(\delta_i = 0)$.

The hazard rate $h(x) := \lim_{\Delta x \rightarrow o}
(1 / \Delta x)\cdot P(x \leq T_i < x + \Delta x | T_i \geq x)$,\linebreak
$x \geq 0$,
has the application-orientated flavour that it can be nicely interpreted
as the instantaneous risk of observing the failure event of interest at
\mbox{time $x$.} In a variety of applications ranging from survival analysis in
a medical context to reliability testing in industrial settings the hazard
rate is thus used extensively as a methodological tool to describe 
variations in risk over time. In these applications qualitative aspects of
the structure of the hazard rate are more interesting than the precise
location and height of peaks of the function. 

The most prominent nonparametric approach to estimate the hazard rate is
given by the kernel estimator with a fixed bandwidth which is defined as
\[
\widehat{h} (x) := \sum_{i=1}^n \frac{\delta_{(i)}}{n-i+1}\cdot 
\frac{1}{b}\cdot K \left( \frac{x-X_{(i)}}{b}\right)
\]
where $\delta_{(i)}$ refers to the censoring indicator corresponding
to the $i$-th element of the order statistic of the observed failure times,
$K(\cdot)$ denotes the kernel function (satisfying standard conditions,
see below), and the bandwidth parameter $b$ has to be positive. Several
variations of this kernel estimator have been suggested
allowing the bandwidth to vary with $x$, for example, the nearest
neighbour kernel estimator (Gefeller and Dette (1992)), the
local-bandwidth kernel estimator (M\"uller and Wang (1994)) or the
variable kernel estimator (Sch\"afer (1985)). For the rest of this
paper attention is focussed on the simplest case of the
fixed-bandwidth kernel estimator. Modifications for the local-bandwidth
and the nearest neighbour kernel estimator are straightforward; however,
the variable kernel estimator poses the additional complexity that
expressions for the $MISE$ have not been derived yet and thus
conditions under which the results of following section might be
transferred to the variable kernel estimator are not yet
explicitly available.

A variety of results on the asymptotic behaviour of $\widehat{h}$ ranging
from different proofs of consistency to sophisticated elaborations on
the optimal order of convergence can be found in the literature. For the
purpose of this paper it is sufficient to restrict the attention to
a result on the asymptotic $MISE$ which can be easily decomposed into
integrals over
a squared bias ($\mu^2(\cdot))$ and a variance part ($\sigma^2 (\cdot)$).
The necessary technical assumptions for the asymptotics to work can be
stated as follows:
\begin{description}
\item[(K)] The kernel function $K(\cdot)$ has to defined on a compact
support $[a, b] \subset {\cal R}$ as a bounded symmetrical probability
density function having a second derivative that is Lipschitz continuous
on $[a, b]$.
\item[(H)] The hazard rate $h$ has to be twice continuously differentiable
and square integrable on ${\cal R}_+$.
\item[(B)] The sequence of bandwidths $b_n$ has to approach zero at a
rate slower than $n^{-1}$. The convergence rate in this situation
is optimised for the bandwidth sequence $b_n = C_0\cdot n^{-1/5}$, with
$C_0 > 0$ denoting a special constant.
\end{description}
Given that (K), (H) and (B) hold it has been shown that for $n\rightarrow
\infty$ the $MISE$ can be written as follows (omitting 
asymptotically vanishing terms):
\begin{eqnarray*}
MISE (\widehat{h}, h)& =& \int E\left[ \widehat{h} (x) - h(x)\right]^2 dx\\
&=&\int \left( 
\left[\frac{b_n^2}{2} h^{\prime \prime}(x) \beta(K)\right]^2 +
\frac{h(x) \alpha(K)}{n b_n (1-F(x)) (1-G(x))}\right)dx\\
&=:& \int \left( \mu^2 (x) + \sigma^2(x)\right) dx
\end{eqnarray*}
where $\alpha(K)$ and $\beta(K)$ denote constants depending on the
kernel function
($\alpha (K) := \int K^2(u)du,\: \beta(K) := \int u^2 K(u) du$).
This result demonstrates the well-known problem in selecting an appropriate
bandwidth, often termed the ``variance-bias trade-off''. For the bias
to decrease one needs to select a small bandwidth, however,
taking this parameter small means automatically an increase in the
variance. The variance-bias trade-off is in accordance with the
intuitive understanding of smoothing as a technique to reveal the underlying
structure of the data by reducing the ``noise'' (variance) at the
expense of some ``oversimplification'' (bias). More details on these
general aspects of smoothing and on the technical details of the derivation
of asymptotic results in a general counting process framework can be
found in the monograph by Andersen et al. (1993).

\section{Asymptotic Properties of Visual Error Criteria Applied to
Hazard Rate Estimation}
The standard asymptotic results on the properties of the kernel estimator
$\widehat{h}$ consider only distances between $\widehat{h}(x)$ and
$h(x)$ at fixed $x \in {\cal R}_+$ and thus exhibit the drawback of
measuring error only vertically as discussed previously. In this section
we analyse the asymptotic behaviour of the quadratic versions of the
visual error criteria that is essentially determined by the asymptotics
of the distances $d((x, \widehat{h}(x)), G_h)$ and $d((x, h(x)), 
G_{\widehat{h}})$, respectively. To this
end, consider for fixed $x_0 \in {\cal R}_+$ the distances as sequences of
nonnegative random variables $D_n^1(x_0)$ and $D_n^2(x_0)$, respectively.
In the context of density estimation Marron and Tsybakov (1995) derived
a result on the convergence in probability of 
$D_n^1(x_0)$ and $D_n^2(x_0)$, respectively, for $n \rightarrow \infty$.
In the proof they used a combination of primarily geometric arguments
that can be directly transferred to the hazard rate context. The only
density-specific step concerned the convergence in probability of the
kernel estimator for the derivative of the density to the true derivative
of the density. The same property, i.e. $\widehat{h}^{\prime} (x_0)
\stackrel{\cal P}{\rightarrow} 
h^{\prime} (x_0)$ for fixed $x_0$ as $n \rightarrow
\infty$, holds for the kernel estimator of the hazard rate as can be checked
easily. Thus, analogous to the density context analysed in Marron and
Tsybakov (1995),
\[
n^{2/5} \cdot \left( D_n^i (x_0) - \frac{| \widehat{h}(x_0) - h(x_0) |}
{\sqrt{1 + \left[ h^{\prime} (x_0)\right]^2}}\right) 
\stackrel{{\cal P}}{\longrightarrow} 0
\]
holds in the hazard rate context for $i=1, 2$ as $n \rightarrow \infty$.
This result allows to build an asymptotic bridge between the expected
squared visual error criteria and the conventional approach as follows:
\begin{eqnarray*}
E\left[ VE_2 (\widehat{h} \rightarrow h)^2\right]& = &
\int d((x, \widehat{h}(x)), G_h)^2 dx\\
&=& \int E\left( D_n^1(x)^2\right) dx\\
&\approx& \int \frac{E\left( [\widehat{h} (x) - h(x)]^2\right)}
{1 + \left[ h^{\prime} (x)\right]^2} dx\\
&\approx& \int \frac{\mu^2 (x) + \sigma^2 (x)}
{1 + \left[ h^{\prime} (x)\right]^2} dx    \quad .
\end{eqnarray*}
Here the first two equalities are given by the definition of the
quantities, the first approximation utilises the asymptotic result on
$D_n^1$ and the last approximation results from plugging in the asymptotic
$MSE$ expression given in the previous section (a similar line of
reasoning leads to the same asymptotic result for 
$E[ VE_2 (\widehat{h} \rightarrow h)^2]$). When defining 
$w (x) := (1 + [h^{\prime} (x)]^2)^{-1}$, the final expression above
can also be viewed as the standard asymptotic representation of a
weighted $MISE$ employing the special weight function $w(\cdot)$.
This shows that expected squared visual error corresponds asymptotically
to a weighted $MISE$, and inspection of the weight function allows to
infer in which situations the two concepts of measuring discrepancy
between curves will contradict each other. For example, for hazard rates
of exponentially distributed failure times
both error concepts give identical answers, but for hazard rates
with regions where $| h^{\prime} (x) |$ is large $E[ VE_2 (\cdot) ^2]$
and $MISE$ can disagree remarkably. This finding corresponds with the
motivating example given in section 2 and the intuitive understanding
of the deficiencies of the $L_2$-approach in situations where
functions reveal several peaks.

Although the asymptotic expressions derived above for $E[VE_2 (\widehat{h} 
\rightarrow h)^2]$ and
$E[VE_2 (h \rightarrow \widehat{h})^2 ]$
are the same, details of the proof in Marron and Tsybakov (1995) suggest that
the asymptotics might come into action later for $VE_2 (h \rightarrow
\widehat{h})$ than for $VE_2 (\widehat{h} \rightarrow h)$. This should
also be reflected in a different finite sample behaviour of the two criteria
which needs further investigation. From the findings on the two
asymmetrical versions of squared visual error it follows immediately that
for the symmetrised criterion $SE_2 (\widehat{h}, h)$ it can
be shown that, for $n \rightarrow \infty$,
$E[SE_2 (\widehat{h}, h)^2]$ tends 
to twice the same weighted $MISE$ as above.

\section{Discussion}
The standard approach when evaluating the performance of a functional 
estimate uses some $L_p$-norm to quantify its distance to the true 
function. As demonstrated in this paper, there are several situations in 
which this approach contradicts the graphical notion of discrepancy 
between curves since all $L_p$-norms consider only vertical distances and 
neglect aspects of qualitative similarities. Thus, the new concept of 
visual error criteria has been discussed as an alternative method to 
evaluate the visual appropriateness of functional estimates. In the 
context of hazard rate estimation from censored data it has been shown 
that application of these criteria corresponds asymptotically to a 
weighted version of the conventional $MISE$. A more detailed analysis of 
the finite properties of visual error criteria is, however, needed to 
expand the knowledge about advantages and disadvantages of this new 
concept. Visual error criteria represent an attractive first step into 
the direction of rethinking the mathematical evaluation of the 
performance of functional estimates, but they require further 
elaboration --- and perhaps some modification --- prior to their routine 
application.

\subsection*{References}

{\small
ANDERSEN, P.K., BORGAN, O., GILL, R.D., KEIDING, N. (1993):
{\em Statistical Models Based on Counting Processes.} Springer-Verlag,
New York.

\medskip
GEFELLER, O., DETTE, H. (1992): Nearest neighbour kernel estimation
of the hazard function from censored data. {\em Journal of Statistical
Computation and Simulation, 43, 93--101.}

\medskip
GEFELLER, O., MICHELS, P. (1992): Nichtparametrische Analyse von
Verweildauern. {\em \"Osterreichische Zeitschrift f\"ur Statistik und
Informatik, 22, 37--59.}

\medskip
H{\"A}RDLE, W. (1991): {\em Applied Nonparametric Regression.}
Cambridge University Press, Boston.

\medskip
HJORT, N.L. (1996): New methods for hazard rate estimation. In:
{\it Proceedings of the 17th Rencontre Franco-Belge de Staticiens,
Universit{\'e} de Marne-la-Vall{\'e}e, p. 20-24.}

\medskip
KOOPERBERG, C., STONE, C.J. (1991): A study of logspline density 
estimation. {\em Computational Statistics \& Data Analysis, 12, 
327--347.}

\medskip
MARRON, J.S., TSYBAKOV, A.B. 
(1995): Visual error criteria for qualitative smoothing.
{\em Journal of the American Statistical Association, 90, 499--507.}

\medskip
M\"ULLER, H.G., WANG, J.L. (1994): Hazard rate estimation under
random censoring with varying kernels and bandwidths.
{\em Biometrics, 50, 61--76.}

\medskip
PARZEN, E. (1962): On the estimation of a probability density and mode.
{\em Annals of Mathematical Statistics, 33, 1065--1076.}

\medskip
ROSENBLATT, M. (1956): Remarks on some nonparametric estimates of a
density function. {\em Annals of Mathematical Statistics, 27, 832--837.}

\medskip
SCH\"AFER, H. (1985): A note on data-adaptive kernel estimation of the
hazard function in the random censorship situation.
{\em Annals of Statistics, 13, 818--820.}

\medskip
WAND, M., JONES, M. (1995): {\em Kernel Smoothing.} Chapman \&
Hall, London.
}
\end{document}